\documentclass[12pt,a4paper]{article}
\usepackage[utf8x]{inputenc}
\usepackage{amsmath}
\usepackage{amsfonts}
\usepackage{amssymb}
\usepackage{amsthm}
\usepackage{mathrsfs}
\usepackage{mathtools}
\usepackage{fancyhdr}
\usepackage{graphicx}
\usepackage[usenames,x11names]{xcolor}
\usepackage{arabtex}
\usepackage{framed}
\usepackage[top=2cm,bottom=2cm,margin=2cm]{geometry}
\usepackage{cancel}
\usepackage{array} 
\usepackage{comment}

\usepackage{hyperref}
\usepackage[numbers,sort&compress]{natbib}

\title{On refinements of two-term Machin-like formulas}
\author{\sc Bakir FARHI \\
National Higher School of Mathematics \\
P.O.Box 75, Mahelma 16093, Sidi Abdellah (Algiers) \\
Algeria \\[1mm]
\href{mailto:bakir.farhi@nhsm.edu.dz}{\tt bakir.farhi@nhsm.edu.dz} \\[1mm]
\url{http://farhi.bakir.free.fr/}
}
\date{}

\let\up=\textsuperscript
\let\epsilon=\varepsilon

\def\Q{{\mathbb Q}}

\def\N{{\mathbb N}}
\def\Z{{\mathbb Z}}

\newcommand{\vabs}[1]{\left\vert{#1}\right\vert}

\theoremstyle{plain}
\numberwithin{equation}{section}
\newtheorem{thm}{Theorem}[section]

\newtheorem{prop}[thm]{Proposition}
\newtheorem{coll}[thm]{Corollary}

\newtheorem{thmn}{Theorem}
 
\newtheorem{propn}[thmn]{Proposition}

\theoremstyle{definition}

\theoremstyle{remark}
\newtheorem{rmk}[thm]{Remark}

\begin{document}
\maketitle

\begin{abstract}
We develop a refinement process for two-term Machin-like formulas: $a_0 \arctan{u_0} + a_1 \arctan{u_1} = \frac{\pi}{4}$ (where $a_0 , a_1 \in \Z$, $u_0 , u_1 \in \Q_+^*$, $u_0 > u_1$) by exploiting the continued fraction expansion of the ratio $\alpha := \frac{\arctan{u_0}}{\arctan{u_1}}$. This construction yields a sequence of derived two-term Machin-like formulas: $a_{- n} \arctan{u_n} + a_{- n + 1} \arctan{u_{n + 1}} = \frac{\pi}{4}$ ($n \in \N$) with positive rational arguments $u_n$ decreasing to zero and corresponding integer coefficients $a_{- n}$. We derive closed forms and estimates for $a_{-n}$ and $u_n$ in terms of the convergents of $\alpha$ and prove that the associated rational sequence $(a_{- n} u_n + a_{- n + 1} u_{n + 1})_n$ converges to $\pi/4$ with geometric decay. The method is illustrated using Euler's two-term Machin-like formula : $\arctan(1/2) + \arctan(1/3) = \pi/4$.
\end{abstract}


\section{Introduction and Notation}\label{sec1}

Throughout this paper, $\N$ denotes the set of positive integers and $\N_0 := \N \cup \{0\}$ the set of non-negative integers. We denote by $(F_n)_{n \in \N_0}$ the usual Fibonacci sequence defined by $F_0 = 0$, $F_1 = 1$, and $F_{n + 2} = F_{n} + F_{n + 1}$ for all $n \in \N_0$. The closed form-expression for $F_n$ ($n \in \N_0$) is given by the Binet formula:
\[
F_n = \dfrac{\Phi^n - \overline{\Phi}^n}{\sqrt{5}} ,
\]
where $\Phi := \frac{1 + \sqrt{5}}{2}$ is the golden ratio and $\overline{\Phi} := - \frac{1}{\phi} = \frac{1 - \sqrt{5}}{2}$ (see e.g., \cite{kos}). Since $\vabs{\overline{\Phi}} < 1$, we have in particular $F_n \sim_{+ \infty} \frac{\Phi^n}{\sqrt{5}}$. Further, we use the Landau big-O notation in its usual sense to describe asymptotic behavior. 

The history of the computation of $\pi$ changed dramatically in 1706 when the English astronomer and mathematician John Machin discovered the famous identity
\begin{equation}\label{eq1-1}
    \frac{\pi}{4} = 4 \arctan\left(\frac{1}{5}\right) - \arctan\left(\frac{1}{239}\right).
\end{equation}
By combining this relation with the Gregory-Leibniz series expansion for the arctangent function,
\begin{equation}\label{eq1-2}
    \arctan x = \sum_{n=0}^{+ \infty} \frac{(-1)^n}{2n+1} x^{2n+1} \qquad (|x| < 1),
\end{equation}
Machin was able to compute the first 100 decimal digits of $\pi$ for the first time in history \cite{Beckmann1971, Tweddle1991}. The validity of Machin's formula, and others of its kind, relies on the fundamental addition and subtraction formulas for the arctangent function:
\begin{equation}\label{eq1-3}
    \begin{split}
        \arctan x + \arctan y & = \arctan\left(\frac{x + y}{1 - x y}\right), \\[1mm]
        \arctan x - \arctan y & = \arctan\left(\frac{x - y}{1 + x y}\right).
    \end{split}
\end{equation}
Consequently, almost all computations of $\pi$ performed between 1700 and 1980 utilized variations of Machin's formula in conjunction with the Taylor series expansion of the arctangent function near zero \cite{Borwein2014, Nishiyama2013}.

In a broader sense, a Machin-like formula is defined as any identity of the form
\begin{equation} \label{eq:general_machin}
    \sum_{i \in I} a_i \arctan u_i = \frac{\pi}{4},
\end{equation}
where $I$ is a finite index set, $a_i \in \mathbb{Z}$, and $u_i \in \mathbb{Q}_+^*$. The computational efficiency of such formulas is directly related to the magnitude of the arguments : as the $u_i$ approach $0$, the convergence of the series (\ref{eq1-2}) accelerates significantly. Consequently, the search for Machin-like formulas involving small arguments remains one of the most effective methods for high-precision calculations of $\pi$. Notable historical examples include the formulas
\begin{align}
    \frac{\pi}{4} & = \arctan\left(\frac{1}{2}\right) + \arctan\left(\frac{1}{3}\right), \label{eq:euler} \\[1mm]
    \frac{\pi}{4} & = 12 \arctan\left(\frac{1}{18}\right) + 8 \arctan\left(\frac{1}{57}\right) - 5 \arctan\left(\frac{1}{239}\right), \label{eq:gauss} \\[1mm]
    \frac{\pi}{4} & = 8 \arctan\left(\frac{1}{10}\right) - \arctan\left(\frac{1}{239}\right) - 4 \arctan\left(\frac{1}{515}\right), \label{eq:simson}
\end{align}
attributed to Euler, Gauss (1863), and Simson (1723), respectively \cite{Tweddle1991,Yamada2025}. The efficacy of this approach was demonstrated in 2002 when Kanada set a world record by computing over one trillion decimal digits of $\pi$ using the following self-checking pair of multi-term formulas \cite{Abrarov2017}:
\begin{align*}
    \frac{\pi}{4} & = 44 \,\arctan\left(\frac{1}{57}\right) + 7 \,\arctan\left(\frac{1}{239}\right) - 12 \,\arctan\left(\frac{1}{682}\right) + 24 \,\arctan\left(\frac{1}{12943}\right) , \\[3mm]
    \frac{\pi}{4} & = 12 \,\arctan\left(\frac{1}{49}\right) + 32 \,\arctan\left(\frac{1}{57}\right) - 5 \,\arctan\left(\frac{1}{239}\right) + 12 \,\arctan\left(\frac{1}{110443}\right) .
\end{align*}
Reflecting the enduring interest in this field, many new identities and methods for generating them have been reported in the modern literature \cite{Abrarov2022, Gasull2023, Yamada2025}.

In this paper, we develop a refinement process for two-term Machin-like formulas
\[
a_0 \arctan{u_0} + a_1 \arctan{u_1} = \dfrac{\pi}{4} ,
\]
with $a_0 , a_1 \in \Z$, $u_0 , u_1 \in \Q_+^*$, and $u_0 > u_1$. By exploiting the continued fraction expansion of the ratio $\frac{\arctan{u_0}}{\arctan{u_1}}$, we construct a sequence of rational positive arguments $u_n$ decreasing to zero together with corresponding integer coefficients $a_{-  n}$ such that
\[
a_{- n} \arctan{u_n} + a_{- n + 1} \arctan{u_{n + 1}} = \dfrac{\pi}{4}
\]
for all $n \in \N_0$. We find closed forms and estimates for $\arctan{u_n}$ and $a_{- n}$ in terms of the convergents of $\frac{\arctan{u_0}}{\arctan{u_1}}$, and prove that the associated rational sequence
\[
\left(a_{- n} u_n + a_{- n + 1} u_{n + 1}\right)_{n \in \N_0}
\]
converges to $\pi/4$ with geometric convergence rate. We illustrate the method by an application for the classical two-term Machin-like formula: $\arctan(1/2) + \arctan(1/3) = \pi/4$.

\section{The results and the proofs}

In what follows, we fix integers $a_0 , a_1$ and positive rational numbers $u_0 , u_1$, with $u_0 > u_1$, such that:
\begin{equation}\label{eq2-1}
a_0 \arctan{u_0} + a_1 \arctan{u_1} = \dfrac{\pi}{4} , 
\end{equation}
and we set 
\[
\alpha := \dfrac{\arctan{u_0}}{\arctan{u_1}} > 1 .
\]
We have the following proposition.

\begin{prop}\label{p1}
The real number $\alpha$ is irrational.
\end{prop}

To prove Proposition \ref{p1}, we need the following well-known result.

\begin{thmn}[{\cite[Corollary 3.12]{niv}}]\label{thmn1}
The only rational values of $\sin(\pi r)$, $\cos(\pi r)$, and $\tan(\pi r)$, for rational $r$, are $0$, $\pm 1/2$, $\pm 1$ for the sine and cosine, and $0$, $\pm 1$ for the tangent.
\end{thmn}

\begin{proof}[Proof or Proposition \ref{p1}]
Assume, for the sake of contradiction, that $\alpha$ is rational. From \eqref{eq2-1}, we obtain:
\begin{equation}\label{eq2-2}
\dfrac{\pi}{4} = a_0 \arctan{u_0} + a_1 \arctan{u_1} = \left(a_0 \dfrac{\arctan{u_0}}{\arctan{u_1}} + a_1\right) \arctan{u_1} = \left(a_0 \alpha + a_1\right) \arctan{u_1} . 
\end{equation}
This implies
\[
\tan\left(\dfrac{1}{4 (a_0 \alpha + a_1)} \pi\right) = u_1 \in \Z .
\]
By Theorem \ref{thmn1}, it follows that 
\[
\tan\left(\frac{1}{4 (a_0 \alpha + a_1)} \pi\right) \in \{0 , 1 , -1\} .
\]
Since $\frac{1}{4 (a_0 \alpha + a_1)} \pi = \arctan{u_1} \in (0 , \frac{\pi}{2})$ (because $u_1 > 0$), we must have 
\[
\frac{1}{4 (a_0 \alpha + a_1)} \pi = \frac{\pi}{4} ,
\]
and therefore $a_0 \alpha + a_1 = 1$. Substituting into \eqref{eq2-2} yields $u_1 = 1$.

By symmetry (interchanging the roles of $u_0$ and $u_1$, and $a_0$ and $a_1$), we similarly obtain $u_0 = 1$, which contradicts the assumption $u_0 > u_1$. Hence $\alpha$ is irrational.
\end{proof}

Now consider the continued fraction expansion of $\alpha$:
\[
\alpha = \left[q_0 ; q_1 , q_2 , \dots\right] ,
\]
where $q_i \in \N$ for all $i \in \N_0$. The infinitude of this expansion follows from the irrationality of $\alpha$ (established in Proposition \ref{p1}). From this continued fraction expansion and the addition and subtraction formulas \eqref{eq1-3} for the arctangent function, there exists a decreasing sequence of positive rational numbers ${(u_n)}_{n \in \N_0}$ such that
\begin{equation}\label{eq-S}
\begin{rcases}
\arctan{u_0} & = q_0 \arctan{u_1} + \arctan{u_2} \\
\arctan{u_1} & = q_1 \arctan{u_2} + \arctan{u_3} \\
~  & \vdots \\
\arctan{u_n} & = q_n \arctan{u_{n + 1}} + \arctan{u_{n + 2}} \\
~ & \vdots
\end{rcases} . \tag{$S$}
\end{equation}
Substituting the first equality of System \eqref{eq-S} into \eqref{eq2-1}, we obtain
\[
\dfrac{\pi}{4} = a_0 \left(q_0 \arctan{u_1} + \arctan{u_2}\right) + a_1 \arctan{u_1} = \left(q_0 a_0 + a_1\right) \arctan{u_1} + a_0 \arctan{u_2} .
\]
Setting
\[
a_{-1} := q_0 a_0 + a_1 \in \Z ,
\]
this becomes
\[
a_{-1} \arctan{u_1} + a_0 \arctan{u_2} = \dfrac{\pi}{4} .
\]
Proceeding similarly by substituting the second equality of \eqref{eq-S}, we obtain
\[
\dfrac{\pi}{4} = a_{-1} \left(q_1 \arctan{u_2} + \arctan{u_3}\right) + a_0 \arctan{u_2} = \left(q_1 a_{-1} + a_0\right) \arctan{u_2} + a_{-1} \arctan{u_3} .
\]
Setting
\[
a_{-2} := q_1 a_{-1} + a_0 \in \Z ,
\]
we get
\[
a_{-2} \arctan{u_2} + a_{-1} \arctan{u_3} = \dfrac{\pi}{4} .
\]
Iterating this process indefinitely leads to the definition of a sequence of integers ${(a_n)}_{n \leq 1}$ given recursively by
\begin{equation}\label{eq2-3}
a_{- n - 1} := q_n a_{- n} + a_{- n + 1} \qquad (\forall n \in \N_0) ,
\end{equation}
which satisfies the general two-term Machin-like formula:
\begin{equation}\label{eq2-4}
a_{- n} \arctan{u_n} + a_{- n + 1} \arctan{u_{n + 1}} = \dfrac{\pi}{4} \qquad (\forall n \in \N_0) 
\end{equation}
(as can be verified by simple induction). As $n$ increases, Formula \eqref{eq2-4} becomes increasingly refined, in the sense that the arguments of the arctangent function strictly decrease. However, these refinements raise two important questions: 
\begin{enumerate}
\item[(1)] Does the sequence ${(u_n)}_n$ converge to $0$?
\item[(2)] Assuming an affirmative answer to Question (1), we have $\arctan{u_n} \sim_{+ \infty} u_n$; so can one deduce from \eqref{eq2-4} that
\[
\lim_{n \to + \infty} \left(a_{- n} u_n + a_{- n + 1} u_{n + 1}\right) = \dfrac{\pi}{4} \, ?
\]
An affirmative answer yields a rational sequence converging to $\pi$, namely
\[
\big\{4 (a_{- n} u_n + a_{- n + 1} u_{n + 1})\big\}_n .
\]
\end{enumerate}
The purpose of this paper is to answer both questions affirmatively. As an application, we consider the case starting from the simplest two-term Machin-like formula: $\arctan(\frac{1}{2}) + \arctan(\frac{1}{3}) = \frac{\pi}{4}$ (where $a_0 = a_1 = 1$, $u_0 = \frac{1}{2}$, and $u_1 = \frac{1}{3}$).

To address Questions (1) and (2) above, let us denote, for each $n \in \N_0$, by $N_n$ and $D_n$ respectively the numerator and the denominator of the positive rational number $[q_0 ; q_1 , q_2 , \dots , q_n]$, so that
\[
\left[q_0 ; q_1 , q_2 , \dots , q_n\right] = \dfrac{N_n}{D_n} \qquad (\forall n \in \N_0) ,
\]
with $N_n , D_n \in \N$ and $\mathrm{gcd}(N_n , D_n) = 1$ for all $n \in \N_0$. According to the elementary theory of continued fractions (see, for example, \cite{khi}), the sequences ${(N_n)}_n$ and ${(D_n)}_n$ are given recursively by
\begin{equation}\label{eq2-5}
\begin{split}
N_{-2} = 0 , N_{- 1} = 1 , \\
D_{-2} = 1 , D_{- 1} = 0 , \\
N_k = q_k N_{k - 1} + N_{k - 2} , \\
D_k = q_k D_{k - 1} + D_{k - 2} 
\end{split} \qquad \begin{matrix}
~~~ \\[8mm]
(\forall k \in \N_0) .
\end{matrix}
\end{equation}

\noindent Moreover, the error incurred when approximating $\alpha$ by its $n$-th convergent $N_n/D_n$ ($n \in \N_0$) satisfies
\begin{equation}\label{eq2-6}
\vabs{\alpha - \dfrac{N_n}{D_n}} \leq \dfrac{1}{D_n D_{n + 1}} .
\end{equation}

The following proposition serves as our main tool for deriving closed forms for $a_{- n}$ and $\arctan{u_n}$ ($n \in \N_0$) in terms of the numbers $N_k$ and $D_k$. Note that this is a special case of a classical result found in \cite{ela}.

\begin{propn}\label{p2}
The set of real sequences ${(x_k)}_{k \geq - 2}$ satisfying the recurrence relation
\[
x_k = q_k x_{k - 1} + x_{k - 2} \qquad (\forall k \in \N_0)
\]
is a two-dimensional real vector space whose basis is given by the pair $\left({(N_k)}_{k \geq - 2} , {(D_k)}_{k \geq - 2}\right)$.
\end{propn}

The closed form of $a_{- n}$ ($n \in \N_0$) in terms of the numbers $N_k$ and $D_k$ is given by the following proposition.
\begin{prop}\label{p3}
For every integer $n \geq - 2$, we have:
\[
a_{- n - 1} = a_0 N_n + a_1 D_n .
\]
\end{prop}

\begin{proof}
Apply Proposition \ref{p2} to the sequence with general term $x_k = a_{- k - 1}$ ($k \geq - 2$). Formula \eqref{eq2-3} ensures that ${(x_k)}_{k \geq -2}$ satisfies the recurrence relation $x_k = q_k x_{k - 1} + x_{k - 2}$ ($\forall k \in \N_0$). Therefore, by Proposition \ref{p2}, $x_n$ ($n \geq -2$) can be expressed as a linear combination with constant coefficients (i.e., independent of $n$) of $N_n$ and $D_n$. Testing the initial values $n = -2$ and $n = -1$ shows that these coefficients are $a_0$ and $a_1$, respectively. This proves the result. 
\end{proof}

From Proposition \ref{p3} we derive the following corollary: 

\begin{coll}\label{c1}
We have
\[
a_{- n - 1} \sim_{+ \infty} \dfrac{\pi/4}{\arctan{u_1}} D_n .
\]
In particular, the integer $a_{- n - 1}$ is positive for all sufficiently large $n$.
\end{coll}

\begin{proof}
For all $n \in \N_0$, Proposition \ref{p3} yields
\[
a_{- n - 1} = a_0 N_n + a_1 D_n = D_n \left(a_0 \dfrac{N_n}{D_n} + a_1\right) .
\]
Since
\[
\lim_{n \to + \infty} \dfrac{N_n}{D_n} = \lim_{n \to + \infty} \left[q_0 ; q_1 , q_2 , \dots , q_n\right] = \left[q_0 ; q_1 , q_2 , \dots\right] = \dfrac{\arctan{u_0}}{\arctan{u_1}}
\]
then
\[
\lim_{n \to + \infty} \left(a_0 \dfrac{N_n}{D_n} + a_1\right) = a_0 \dfrac{\arctan{u_0}}{\arctan{u_1}} + a_1 = \dfrac{a_0 \arctan{u_0} + a_1 \arctan{u_1}}{\arctan{u_1}} = \dfrac{\pi/4}{\arctan{u_1}} ,
\]
by \eqref{eq2-1}. The asymptotic formula follows.
\end{proof}

We now determine a closed form for $\arctan{u_n}$ ($n \in \N_0$) in terms of the numbers $N_k$ and $D_k$. We have the following proposition.  

\begin{prop}\label{p4}
For all $n \in \N_0$, we have
\[
\arctan{u_n} = (-1)^n \left(D_{n - 2} \arctan{u_0} - N_{n - 2} \arctan{u_1}\right) .
\]
\end{prop}

\begin{proof}
Consider temporarily the sequence ${(x_k)}_{k \geq -2}$ defined by:
\[
x_k := (-1)^k \arctan{u_{k + 2}} \qquad (\forall k \geq - 2) .
\]
The equations of System \eqref{eq-S} imply that the sequence ${(x_k)}_{k \geq -2}$ satisfies the recurrence relation:
\[
x_k = q_k x_{k - 1} + x_{k - 2} \qquad (\forall k \in \N_0) .
\]
By Proposition \ref{p2}, $x_n$ ($n \in \N_0$) can therefore be expressed as a linear combination with constant coefficients (i.e., independent of $n$) of $N_n$ and $D_n$. Evaluating at $n = -2$ and $n = -1$ shows that these coefficients are $(- \arctan{u_1})$ and $\arctan{u_0}$, respectively. Hence
\[
x_n = \left(- \arctan{u_1}\right) N_n + \left(\arctan{u_0}\right) D_n \qquad (\forall n \in \N_0) ,
\]
which yields the desired formula.
\end{proof}

An important consequence of this last proposition is the following corollary, which in particular provides an affirmative answer to Question (1).

\begin{coll}\label{c2}
For all $n \in \N$, we have
\[
\arctan{u_n} \leq \dfrac{\arctan{u_1}}{D_{n - 1}} .
\]
In particular, 
\[
\lim_{n \to + \infty} u_n = 0 .
\]
\end{coll}

\begin{proof}
The inequality of the corollary is trivially true for $n = 1$. Let $n \geq 2$ be an integer. By Proposition \ref{p4},
\[
\arctan{u_n} = (-1)^n \left(D_{n - 2} \arctan{u_0} - N_{n - 2} \arctan{u_1}\right) = (-1)^n D_{n - 2} \left(\arctan{u_1}\right) \left(\alpha - \dfrac{N_{n - 2}}{D_{n - 2}}\right) . 
\]
Using Estimate \eqref{eq2-6}, we obtain
\begin{align*}
\arctan{u_n} = \vabs{\arctan{u_n}} = D_{n - 2} \left(\arctan{u_1}\right) \vabs{\alpha - \dfrac{N_{n - 2}}{D_{n - 2}}} & \leq D_{n - 2} \left(\arctan{u_1}\right) \cdot \dfrac{1}{D_{n - 2} D_{n - 1}} \\
& = \dfrac{\arctan{u_1}}{D_{n - 1}} ,
\end{align*}
which proves the required inequality. As $n \to + \infty$, since $D_{n - 1} \to + \infty$, we conclude that $\arctan{u_n} \to 0$, and hence $u_n \to 0$. This completes the proof.
\end{proof}

We now turn to the following corollary, which provides a precise affirmative answer to Question (2).

\begin{coll}\label{c3}
As $n \to + \infty$, we have
\[
a_{- n} u_n + a_{- n + 1} u_{n + 1} = \dfrac{\pi}{4} + O\left(\dfrac{1}{D_{n - 1}^2}\right) .
\]
In particular,
\[
\lim_{n \to + \infty} \left(a_{- n} u_n + a_{- n + 1} u_{n + 1}\right) = \dfrac{\pi}{4} .
\]
\end{coll}

\begin{proof}
As $n \to + \infty$, since $u_n \to 0$ (by Corollary \ref{c2}), we may use the second-order Taylor expansion of the tangent function near $0$, yielding
\[
u_n = \tan\left(\arctan{u_n}\right) = \arctan{u_n} + O\left((\arctan{u_n})^3\right) .
\]
By Corollary \ref{c2}, this gives
\begin{equation}\label{eq2-7}
u_n = \arctan{u_n} + O\left(\dfrac{1}{D_{n - 1}^3}\right) .
\end{equation}
On the other hand, Corollary \ref{c1} implies
\begin{equation}\label{eq2-8}
a_{- n} = O\left(D_{n - 1}\right) .
\end{equation}
Combining \eqref{eq2-7} and \eqref{eq2-8}, we obtain
\begin{align*}
a_{- n} u_n & = a_{- n} \arctan{u_n} + O\left(\dfrac{1}{D_{n - 1}^2}\right) \\[-4mm]
\intertext{and similarly,} \\[-12mm]
a_{- n + 1} u_{n + 1} & = a_{- n + 1} \arctan{u_{n + 1}} + O\left(\dfrac{D_{n - 2}}{D_n^3}\right) = a_{- n + 1} \arctan{u_{n + 1}} + O\left(\dfrac{1}{D_{n - 1}^2}\right) .
\end{align*}
Therefore,
\begin{align*}
a_{- n} u_n + a_{- n + 1} u_{n + 1} & = a_{- n} \arctan{u_n} + a_{- n + 1} \arctan{u_{n + 1}} + O\left(\dfrac{1}{D_{n - 1}^2}\right) \\
& = \dfrac{\pi}{4} + O\left(\dfrac{1}{D_{n - 1}^2}\right) \qquad (\text{according to \eqref{eq2-4}}) ,
\end{align*}
which confirms the first assertion of the corollary. The second part follows from the trivial fact that $D_{n - 1} \to + \infty$ as $n \to + \infty$.
\end{proof}

\begin{rmk}
Since all partial quotients $q_k$ ($k \in \N_0$) are positive integers, we have for all $k \in \N_0$:
\[
D_k = q_k D_{k - 1} + D_{k - 2} \geq D_{k - 1} + D_{k - 2} .
\]
A simple induction then shows that
\[
D_{n - 1} \geq F_n \qquad (\forall n \in \N_0) .
\]
Since $F_n \sim_{+ \infty} \frac{\Phi^n}{\sqrt{5}}$ (as recalled in Section \ref{sec1}), it follows that
\[
\dfrac{1}{D_{n - 1}^2} = O\left(\dfrac{1}{F_n^2}\right) = O\left(\dfrac{1}{\Phi^{2 n}}\right) \qquad (\text{as } n \to + \infty) .
\]
Thus, the estimate in Corollary \ref{c3} shows that the rational sequence
\[
{\left(a_{- n} u_n + a_{- n + 1} u_{n + 1}\right)}_n
\]
converges to $\pi / 4$ at least as fast as a geometric sequence with ratio $1 / \Phi^2 \simeq 0.382$.
\end{rmk}

\subsection*{Concrete application.} We conclude by applying our refinement algorithm to the classical two-term Machin-like formula
\begin{equation}\label{eq2-9}
\arctan\left(\frac{1}{2}\right) + \arctan\left(\frac{1}{3}\right) = \dfrac{\pi}{4} .
\end{equation}
Here $a_0 = a_1 = 1$, $u_0 = \frac{1}{2}$, and $u_1 = \frac{1}{3}$. We obtain
\[
\alpha = \dfrac{\arctan(1/2)}{\arctan(1/3)} = \left[1 ; 2 , 3 , 1 , 2 , 1 , 4 , \dots\right] .
\]
Thus, the sequence ${(q_n)}_{n \in \N_0}$ begins with
\[
q_0 = 1 ~,~ q_1 = 2 ~,~ q_2 = 3 ~,~ q_3 = 1 ~,~ q_4 = 2 ~,~ q_5 = 1 ~,~ q_6 = 4 ~,~ \dots .
\]
The terms of the sequence ${(u_n)}_n$ (from $n = 2$) are obtained recursively using System \eqref{eq-S}:
\[
u_2 = \dfrac{1}{7} ~,~ u_3 = \dfrac{3}{79} ~,~ u_4 = \dfrac{24478}{873121} ~,~ \dots .
\]
The sequence ${(a_n)}_n$ for $n < 0$, obtained from \eqref{eq2-3}, begins as
\[
a_{-1} = 2 ~,~ a_{-2} = 5 ~,~ a_{-3} = 17 ~,~ a_{-4} = 22 ~,~ a_{-5} = 61 ~,~ a_{-6} = 83 ~,~ \dots .
\]
The first few refined Machin-like formulas \eqref{eq2-4} are therefore
\begin{align*}
2 \arctan\left(\dfrac{1}{3}\right) + \arctan\left(\dfrac{1}{7}\right) & = \dfrac{\pi}{4} , \\
5 \arctan\left(\dfrac{1}{7}\right) + 2 \arctan\left(\dfrac{3}{79}\right) & = \dfrac{\pi}{4} , \\
17 \arctan\left(\dfrac{3}{79}\right) + 5 \arctan\left(\dfrac{24478}{873121}\right) & = \dfrac{\pi}{4} ,
\end{align*}
and so on. Finally, the corresponding rational approximations of $\pi$, 
\[
r_n := 4 \left(a_{- n} u_n + a_{- n + 1} u_{n + 1}\right) \qquad (n \in \N_0) ,
\]
are given by
\begin{align*}
r_0 & = \dfrac{10}{3} = 3.333\dots , \\[2mm]
r_1 & = \dfrac{68}{21} = 3.238\dots , \\[2mm]
r_2 & = \dfrac{1748}{553} = 3.1609\dots , \\[2mm]
r_3 & = \dfrac{216791924}{68976559} = 3.1429\dots ,~ \text{etc} . 
\end{align*}

\rhead

\end{document}